\documentclass[preprint,12pt]{elsarticle}

\graphicspath{ {./figures/} }
\usepackage{subfig}
\usepackage{xcolor}
\usepackage{amssymb}
\usepackage{amsmath}

\begin{document}

\begin{frontmatter}

\title{Emergence of Order in Dynamical Phases in Coupled Fractional Gauss Map}

\author[1]{Sumit S. Pakhare }
\ead{sumitspakhare@gmail.com}

\author[2]{Varsha Daftardar-Gejji
}
\ead{vsgejji@gmail.com}

\author[1]{Dilip S. Badwaik }
\ead{badwaik_ds@rediffmail.com}

\author[3]{Amey Deshpande }
\ead{2009asdeshpande@gmail.com}

\author[4]{Prashant M. Gade\corref{cor1}} 
\corref{cor1} 
\ead{prashant.m.gade@gmail.com}
\cortext[cor1]{Corresponding Author}

\address [1]{Department of Physics, Kamala Nehru Mahavidyalaya, Nagpur}
\address [2]{Department of Mathematics, Savitribai Phule Pune University, Pune}
\address [3]{Department of Mathematics, Indian Institute of Technology, Bombay}
\address [4]{Department of Physics, RTM Nagpur University, Nagpur}
\begin{abstract}
%% Text of abstract
Dynamical behaviour of discrete dynamical systems has been
investigated extensively in the past few decades.
However, in several
applications, long term memory plays an important role in the evolution of dynamical variables. 
The definition of discrete maps has recently been extended to fractional maps to model such situations. 
We extend this definition to a spatiotemporal system.
We define a coupled map lattice
on different topologies, namely,
one-dimensional coupled map lattice, globally coupled system
and small-world network.
The spatiotemporal patterns in the fractional system are more ordered. 
In particular, synchronization is observed over a large parameter region.
For integer order coupled map lattice in one dimension, synchronized
periodic states with a period greater than one are not obtained.
However, we observe synchronized periodic states with 
period-3 or period-6 in one dimensional coupled
fractional maps even for a large lattice. 
With nonlocal coupling, the synchronization is reached 
over a larger
parameter regime.
In all these cases, the standard deviation decays as power-law in time
with the power same as fractional-order.
The physical significance of such studies is also discussed.
\end{abstract}


\end{frontmatter}

%% \linenumbers

%% main text
\section{Introduction}
\label{}
Discrete fractional calculus has a long history starting with Leibniz (1695) 
and is touched upon by almost every eminent mathematician. The recent revival of
interest in fractional calculus is due to the successful application of 
fractional calculus in different fields such as quantum mechanics, 
electromagnetics, bioengineering, signal processing and many more
\cite{atici2009initial,sheng2011fractional,oldham1974fractional,samko1993fractional,voyiadjis2019brain,podlubny1999fractional}.
Both low-dimensional, as well as high dimensional spatiotemporal
systems, are studied in this context.  Fractional partial differential
equations are also proposed\cite{PhysRevE.66.056108,doi:10.1142/10541,BHRAWY2017197}. 
In this paper, we propose a system of coupled fractional Gauss maps and study the phase transitions in it. 
In coupled map lattice, we study
finite-dimensional difference equations coupled by Laplace operator. 
We generalize
the difference equations to fractional-order.

Several systems in nature have inherent long term memory in it \cite{sumelka2017hyperelastic,nauk1983bulletin}. 
Rheological systems, geophysical systems or fracture dynamics have been 
modelled by 
integrodifferential equations for a long time 
\cite{magin2006fractional,daftardar2013fractional,atici2010modeling,
rudolf2000applications,klimek2001fractional,drapaca2012fractional}.
Fractional differential equations 
offer an alternative modeling scheme
\cite{yang2012advanced,metzler2003fractional}. 
By now, we know several cases where fractional differential equations have been successful in providing a qualitative and quantitative understanding of the system\cite{tarasov2011fractional,lazopoulos2006non}. 

Fractional calculus has been
successfully applied to spatiotemporal systems as well. 
Anomalous random walks have been observed experimentally
since 1921 and several theories have been proposed to understand those. One
of  the recent models uses fractional diffusion equations to model it\cite{C4CP03465A,
METZLER20001,Wu2015}. The fractional diffusion equations 
have also been useful in modelling the behaviour of random walkers in an expanding 
medium\cite{yuste2016diffusion}.

In compressional and shear waves in sediments,
a near-linear variation of attenuation with frequency is
observed over some range. The standard Biot theory predicts
attenuation that increases and frequency squared or either levelling 
off or increasing as square root at high frequencies.
An alternative theory explaining thus
attenuation was proposed by Holm and Pandey. It has time-domain
memory operator equivalent to a fractional derivative operator.
\cite{holm2016wave}. Contribution of ions to electrical impedance in
an electrolytic cell is related to the anomalous diffusion process and
memory effects.  Evengelista {\it{et al}} proposed analytical 
solutions of  fractional diffusion equations  to explain it
\cite{evangelista2011anomalous}. Decomposition of
supersaturated solid solution can lead to  the formation of clusters and
precipitates  of atoms and
defects. They can change many material properties in a significant manner. 
Their kinetics is
described by the diffusion-limited process. Sibatov and Svetukhin
developed generalized equations of subdiffusion-limited
growth using
fractional derivatives to model this situation
\cite{sibatov2015fractional}. Another example 
where fractional calculus has helped to model
physical phenomena is heat conduction. Classical Fourier
law and heat conduction equation are inadequate for describing heat
conduction in several materials.
Time-nonlocal generalization using
fractional calculus was 
proposed by Povstenko in this context\cite{povstenko2016fractional}. 
Fractional calculus has also been used for spatiotemporal systems
such as reaction-diffusion systems\cite{owolabi2016mathematical}. 
Thus fractional equations have ben useful in modelling several
physical systems including spatiotemporal systems.

We can construct a spatiotemporal system which comprises of units which follow
fractional dynamics and such system can shed light on dynamics in 
several naturally occurring spatiotemporal systems.

A system consisting of two coupled fractional Henon maps was proposed in \cite{Liu2016} and the synchronization was studied. 
Though coupled fractional logistic map
has been studied in \cite{zhang2017spatiotemporal}, the formulation is
the same as usual coupled map lattice. The only difference is that
a variant of the logistic map is used as an on-site map 
(see eq. 5 of \cite{zhang2017spatiotemporal}). 
Thus the entire theory of
coupled map lattice can be used in this case.
Our formulation has a long-term memory which
is absent in standard definition of coupled map lattice.
We observe that this helps in establishing
long-range order even in low dimensional systems.
Our definition reduces to a single fractional map in the absence of coupling. We 
study  this formulation
for a few topologies such as one-dimensional lattice, global coupling
and small-world lattice. The prescription is generic enough to be extended
to any underlying topology.

While fractional differential equations have been 
extremely successful in studies of several physical systems,
their numerical simulation is tedious and time-consuming.
It is easier to simulate fractional difference equations.
We may be able to spot the essential characteristics of
the system using such modelling. 
We expect a generic nature of the dynamical evolution of systems with long-term memory that can be learned using these discrete-time models. 
There have been very few investigations in such systems or their spatially extended counterpart.  We study such systems
with local as well as nonlocal coupling. Such modelling may
capture some phenomena in continuous time spatiotemporal systems with
fractional-order.

From the engineering viewpoint,
we note that dynamical systems with long-term memory have been studied in the context
 of control
of chaos. Socolor and co-workers \cite{socolar1994stabilizing} studied systems with
 extended time-delayed
feedback. In this system, the evolution has a component from the feedback signal
 given by 
$\left(1-R\right)\sum_{j=1}^{t} R^{\textit{k}-1} \textit{x}\left(t-\textit{k}\right)$.
Such feedback has been found extremely useful in controlling certain models as well as experimental systems. 
In this case, the simulation can be simplified by introducing another variable and
converting it to an 
effectively two-dimensional system since the kernel is exponential.
In our case, the kernel is a power-law and such simplification is not available. 
However, there have been advances in simulating electrical circuits of non-integer
order and it is possible that such kernels can be experimentally
realized \cite{gomez2018fundamental}.

We follow the definition in \cite{gejji2016}. Deshpande and Daftardar-Gejji
define a fractional difference operator and define the evolution 
as, 
\begin{equation}
x(t)=x\left(0\right)+\frac{1}{\varGamma \left( \alpha \right) }\sum_{j=1}^{t} 
\frac{\varGamma \left(t-j+\alpha \right)}{\varGamma \left(t-j+1 \right)} 
\left[f\left(x\left(j-1\right)\right)-x\left(j-1\right)\right], \label{eq1}
\end{equation}
where, t is an integer, $x\left(0\right)$ is an initial condition, 
$0<\alpha<1$ and $f\left(x\right)$ is a difference equation 
or map. To simplify this notation, we define,
\begin{equation}
g_{\alpha}\left(t-j\right)= \frac{\varGamma \left(t-j+\alpha \right)}
{\varGamma \left(t-j+1 \right)}. \label{eq2}
\end{equation}
Thus
\begin{equation}
x(t)=x\left(0\right)+\frac{1}{\varGamma \left( \alpha \right) }\sum_{j=1}^{t} 
g_{\alpha}(t-j) 
\left[f\left(x\left(j-1\right)\right)-x\left(j-1\right)\right], \label{eq1}
\end{equation}

Often-employed definition of one-dimensional coupled map lattice is following.
Let $x\left(i,t\right)$ be a real variable associated with site $i$ at time $t$. 
The evolution is usually given by,
\begin{equation}
x\left(i,t\right)=\left(1-\epsilon \right)f\left(x\left(i,t-1\right)\right)+\frac{\epsilon }
{2}\left[ f\left(x\left(i+1,t-1\right)\right)+f\left(x\left(i-1,t-1\right)\right)\right], \label{eq3}
\end{equation}
where, $1<i<N$ and $\epsilon$ is the coupling strength.

We generalize this definition 
and define one-dimensional coupled fractional maps in the
next section. The uncoupled system evolves like
$N$ independent fractional maps, {\it{i.e.}} for $\epsilon=0$, the 
evolution at each site is described by a single fractional map.
We also study synchronization and other
 spatiotemporal patterns in this one-dimensional 
 locally coupled fractional maps.
We extend this definition to globally coupled maps in the third section. 
Similar results are obtained for small-world networks.
For non-local couplings, 
we observe synchronization over a large range of parameter space. 
The approach to synchronization is not exponential but a power-law. 
The power is found to be related to the fractional-order parameter. 
With local coupling, synchronized periodic states are observed.
There are clear differences even in the qualitative behaviour of dynamics
in presence of long-term memory.

\section{Coupled fractional Gauss maps}
The discrete Gauss map $f(x)$ is defined as
\begin{equation}
f\left(x\right)=\exp \left(-\nu x^{2}\right)+\beta. \label{eq7}
\end{equation}
We fix $\nu = -7.5$ and vary $\beta$.
The discrete fractional Gauss map according to the definition 
in \cite{gejji2016} is given as,
\begin{equation}
x(t)=x\left(0\right)+\frac{1}{\varGamma \left( \alpha \right) }\sum_{j=1}^{t} 
g_{\alpha}(t-j)\left[\exp 
\left(-7.5 \left(x\left(j-1\right)^{2}\right)\right)+\beta -x\left(j-1\right)\right]. \label{eq9}
\end{equation}

Few comments are in order. Popular maps such as logistic maps or tent maps 
are a function
of unit interval onto itself and coupled map lattice is defined in a way such that coupling keeps the range of values within the unit interval [0,1]. 
However, the fractional variant of this map does not keep the range of
values in the same interval \cite{wu2014discrete}. For some parameter values and initial
conditions, the system blows up and values tend to infinity. We choose Gauss map since 
$f\left(x\right)$ is bounded for any $x$ and is thus more stable. 
Apart from the Gauss map, Daftardar-Geiji and Deshpande studied Bernoulli map in their work. 

We assume periodic boundary conditions and
define the evolution of coupled fractional maps as,
\begin{equation}
x\left(i,t\right)=x\left(i,0\right)+\frac{1}{\varGamma \left(\alpha\right)}
\sum_{j=1}^{t} g_{\alpha}
\left(t-j\right)
G_1(
x\left(i,j-1\right),
x\left(i+1,j-1\right),
x\left(i-1,j-1\right)
) 
\label{eq4}
\end{equation}
where,
\begin{equation}
G_1(a,b,c)
=\left(1-\epsilon \right)f(a)+\frac{\epsilon}{2}\left[f(b)+f(c) \right] -a. \label{eq6}
\end{equation}

This definition matches with the extension of the fractional map to a two-dimensional 
Henon map \cite{podlubny1999fractional}.
For $\epsilon =0$, 
the bifurcation diagram of a single fractional map is reproduced.
This is a high dimensional system (even for a single map) and the attractor may 
depend on initial conditions.
We take random initial conditions. We find that the generic nature of attractor
does not change much with initial conditions as long as they are random. 
We note that the computations become much faster if we store
values of $g_{\alpha}\left( m\right)$ for a given value of $\alpha$ for $m\le T$ 
to simulate the system for $T$ time-steps. These values of ratios of
gamma functions were computed in high precision using Mathematica.

The above definition allows the possibility for the existence of a synchronized state. 
If $x(i,0)=c$ for all $i$,
the values will continue to be synchronized for all times. However, strictly
speaking, synchronized fixed point or synchronized periodic state is not an absorbing state
(unlike usual coupled map lattice).
The future depends on all past values. So even if the values get synchronized at some
time $T$, they need not stay synchronized since past values are different.
Nonetheless, we observe spatial synchronization over a reasonable parameter
range because the weight of
past values reduces as a power-law. 
Asymptotically, we observe a state which is synchronized to an arbitrary
precision.

We present results for three values of $\alpha$ {\it{i.e.}}, 
 $0.4$, $0.6$ and $0.8$. 
We study dynamics for various values of   $\epsilon$  and $\beta$.
For these maps,
$\alpha$ is the fractional-order and $\epsilon$ is the coupling strength. 
We vary $\beta$ between $-0.40$ to $-0.50$ 
which is the critical range. We consider $N=100$
and simulate the model for
$5\times 10^4$ time steps. We use heat maps to
observe and demonstrate spatiotemporal patterns in this system.
The state of all the maps at every time step can be plotted with 
the help of heat maps. In addition to heat maps, we study spatial
snapshot of the profile, time series of values of a given site, 
time-series of mean-field and time series of standard
deviation of mean-field to learn more about dynamics.

For $\alpha = 0.4$, and  $\beta=-0.40$, 
all the maps stay unsynchronized over the observed time. But as $\beta$ 
approaches $-0.50$ the domains of almost synchronized values become bigger and bigger. 
These domains come back to their values after period-3. For  
$N=100$, we observe an almost 3-periodic state after roughly
$ 2 \times 10^4$  time steps.
The maps change their state slightly. However, 
the decrease in the standard deviation of mean is a clear
indicator of approach to synchronization.
In Fig. 1, we show heat map for  $\beta=-0.40$ and $\beta=-0.50$.
It is clear as we have spatiotemporal chaos in one case while
the system reaches perfect order in the second case.
For visual clarity, we have plotted the heat map with only those
values of $x(i,t)$ for which $\mod(t,3)=0$.
The reason is that the underlying periodicity is 3.

Same periodicity and synchronization is observed for 
$\alpha=0.6$ as well. 
We plot heat maps of the system for two values of $\beta$, namely
$\beta=-0.40$ and $\beta=-0.44$ in Fig. 2.
Again, there is a clear transition from spatiotemporal chaos in
one case and a synchronized state in the other.
In this case, also, we have plotted heat map
where the $x(i,t)$ is plotted only for $\mod(t,3)=0$.

For $\alpha=0.8$, the periodicity of the synchronized state
changes to 6.
In Fig. 3, we have shown heat map at 
 $\beta=-0.40$ while for $\beta=-0.45$ 
where the $x(i,t)$ is plotted only for $\mod(t,6)=0$.
We observe a clear transition from spatiotemporal chaos in one case
to a synchronized period-6 state for $\alpha=0.8$.
As the number of maps increases, time taken for synchronization increases.
In Fig. 4, we have shown heat maps for systems with $N=200$
for same parameter values.

The question is whether the synchronization can occur in 
the thermodynamic limit.
Though we do not observe exact synchronization for period-6 states for large
lattices, we believe that synchronization may occur for
the period-3 state even in the thermodynamic limit. 
We investigate the case $\alpha=0.6,  \beta=-0.45$ and $\alpha
=0.4, \beta=-0.55$ in
detail. We compute 
time $T_N$ required for coupled maps to reach a state in which 
$\max(x(i,T_N))-\min(x(i,T_N)) <0.01$. We average over 
20-100 configurations. The average time required $T_N$ 
is plotted as a function of system size $N$ and it shows a 
power-law with exponent 1.22. 
For $\alpha=0.4$, the synchronization time scales with
system size as a power-law again and the
exponent is 0.96. It is shown in Fig. 5.
(For $\alpha=0.8$, the dynamics is multistable, some configurations
reach synchronization and some do not reach synchronization. 
Hence we have not studied the case  $\alpha=0.8$.)
It is clear
that the synchronization can be achieved even in the thermodynamic limit
for a period-3 state for $\alpha=0.4$ and $\alpha=0.6$.
This is quite a unique behaviour. We have not seen reports of synchronized periodicity
even for period-2 in coupled map lattice even for lattice size of 50. 
The reason is the following:
if there is a 
periodicity, different sites of the CML settle down at different phases of this periodicity. 
Thus the spatial synchronization is not obtained with local connections alone 
except for synchronized fixed point
\cite{PhysRevLett.58.2155}. 

In this system, the entire memory of all previous time-steps is involved in the computation. The consequence is that
simulation time scales nonlinearly 
with time-steps.
Simulating very large lattices for a very long time is not practical
for this system.

\begin{figure}%{.3\linewidth}
\centering
\subfloat[]{{\includegraphics[width=\linewidth]{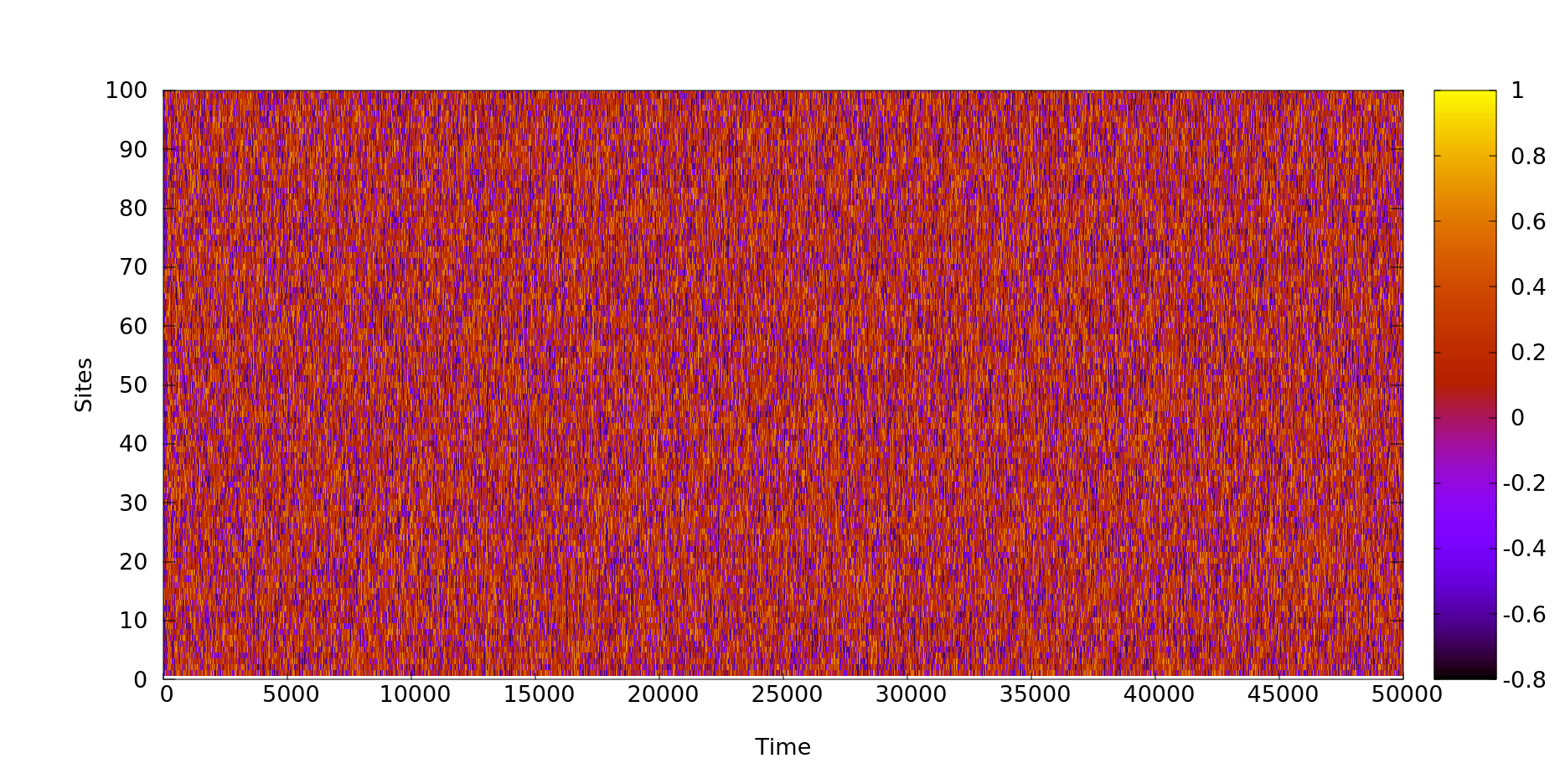} }}%
\qquad
\subfloat[]{{\includegraphics[width=\linewidth]{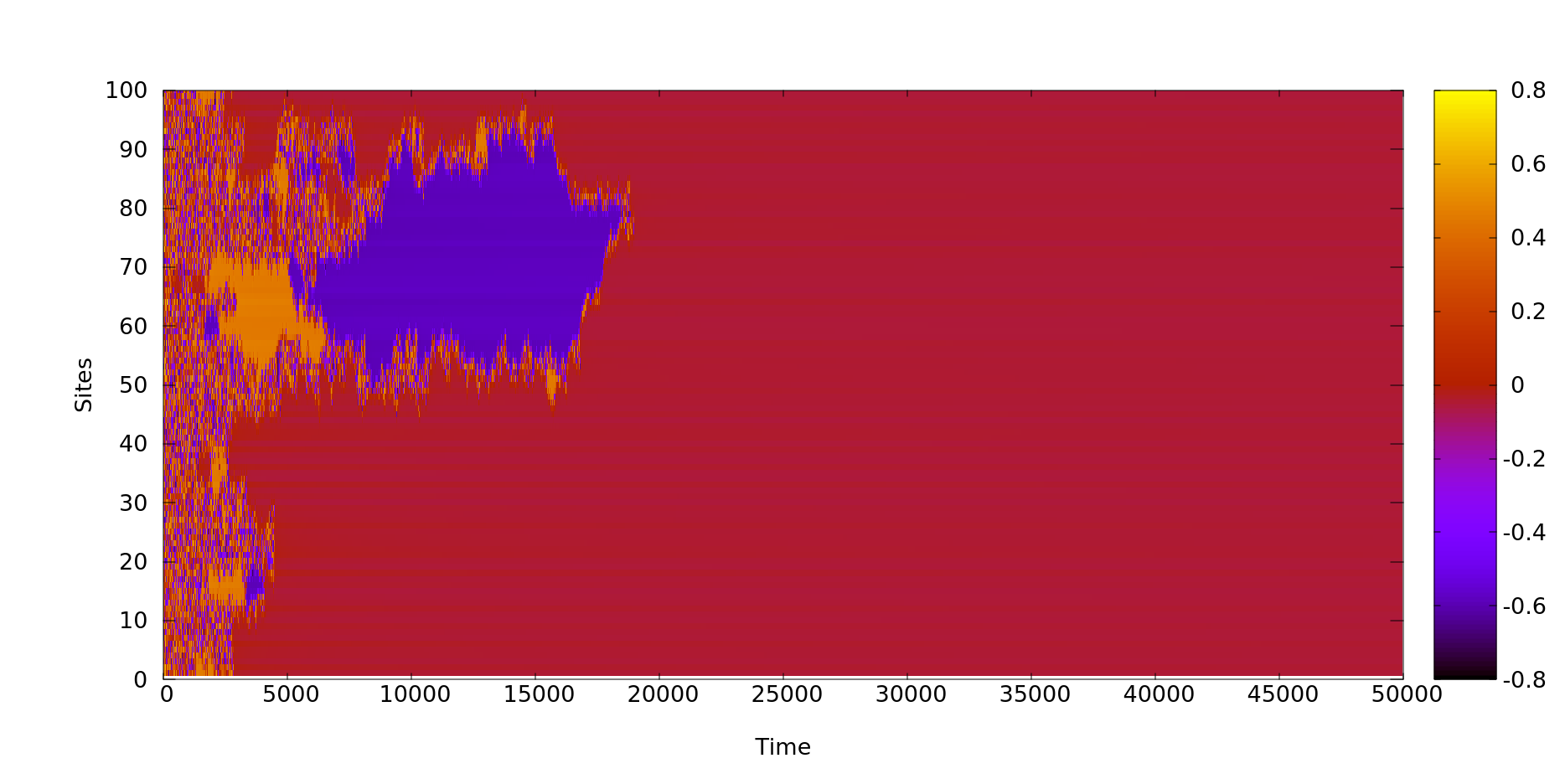} }}%
\caption{Heat map demonstrating spatiotemporal evolution of 
 the system for (a) $\alpha=0.4$ and $\beta=-0.40$
showing unsynchronized behavior for all time-steps while for (b) $\alpha=0.4$ and 
$\beta=-0.50$ synchronization 
after 20000 time-steps is observed.}%
\label{fig:1}%
\end{figure}

\begin{figure}%{.3\linewidth}
\centering
\subfloat[]{{\includegraphics[width=\linewidth]{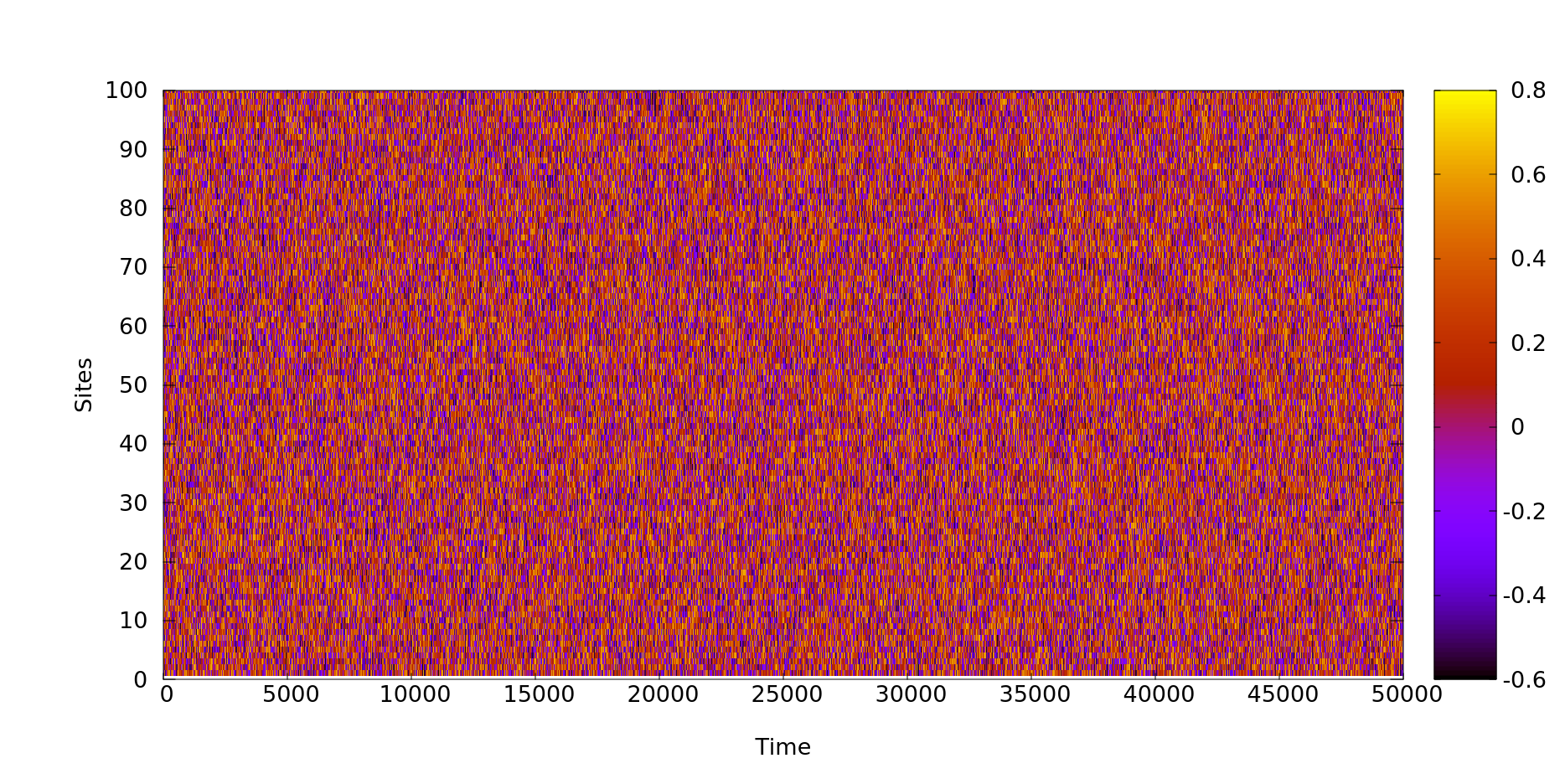} }}%
\qquad
\subfloat[]{{\includegraphics[width=\linewidth]{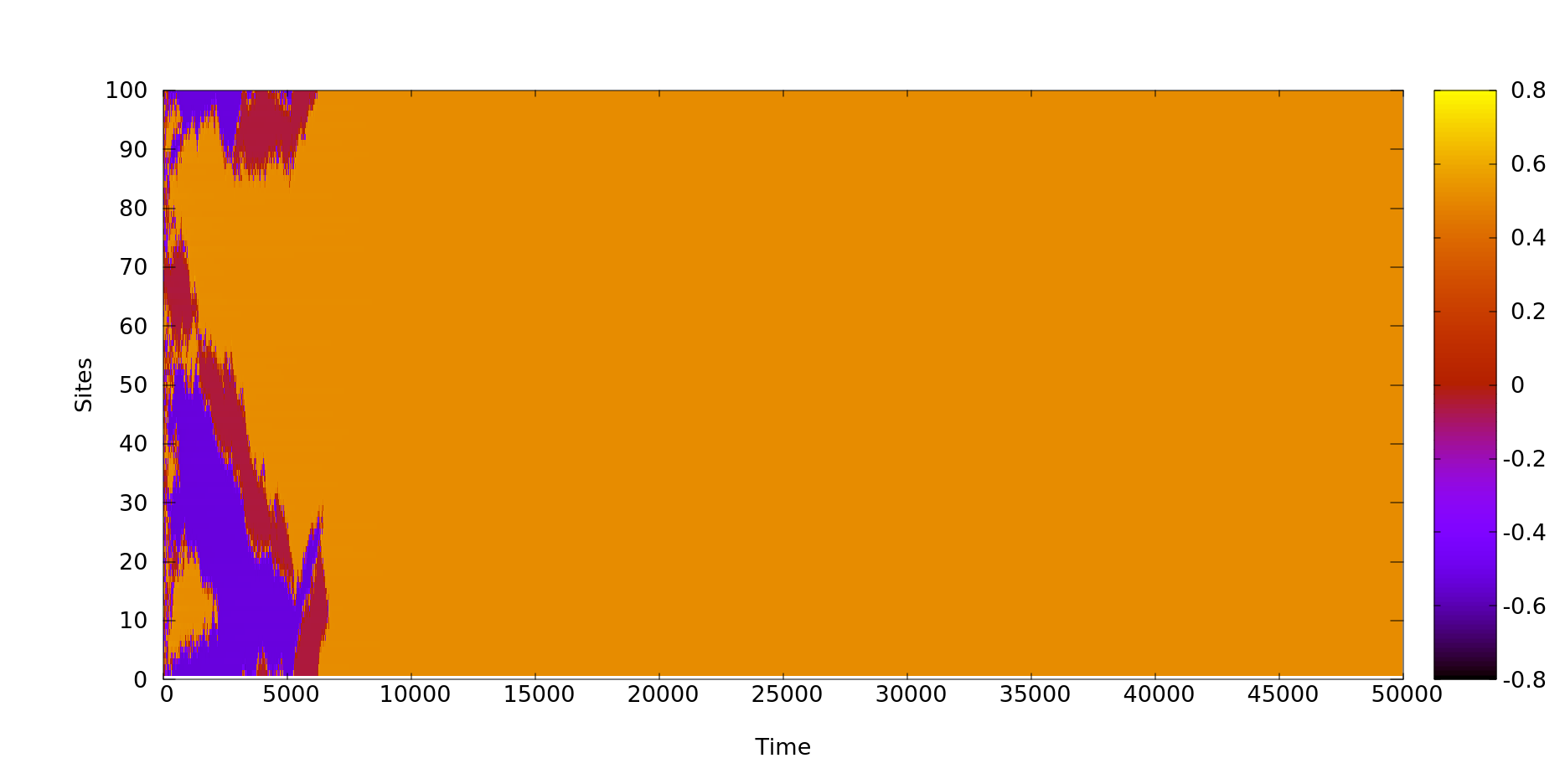} }}%
\caption{Heat map demonstrating spatiotemporal evolution of the system
for (a) $\alpha=0.6$ and $\beta=-0.40$, showing 
unsynchronized behavior for all time-steps while for (b) $\alpha=0.6$ and $\beta=-0.44$
system synchronizes after first 7000 time steps.}%
\label{fig:2}%
\end{figure}

\begin{figure}%{.3\linewidth}
\centering
\subfloat[]{{\includegraphics[width=\linewidth]{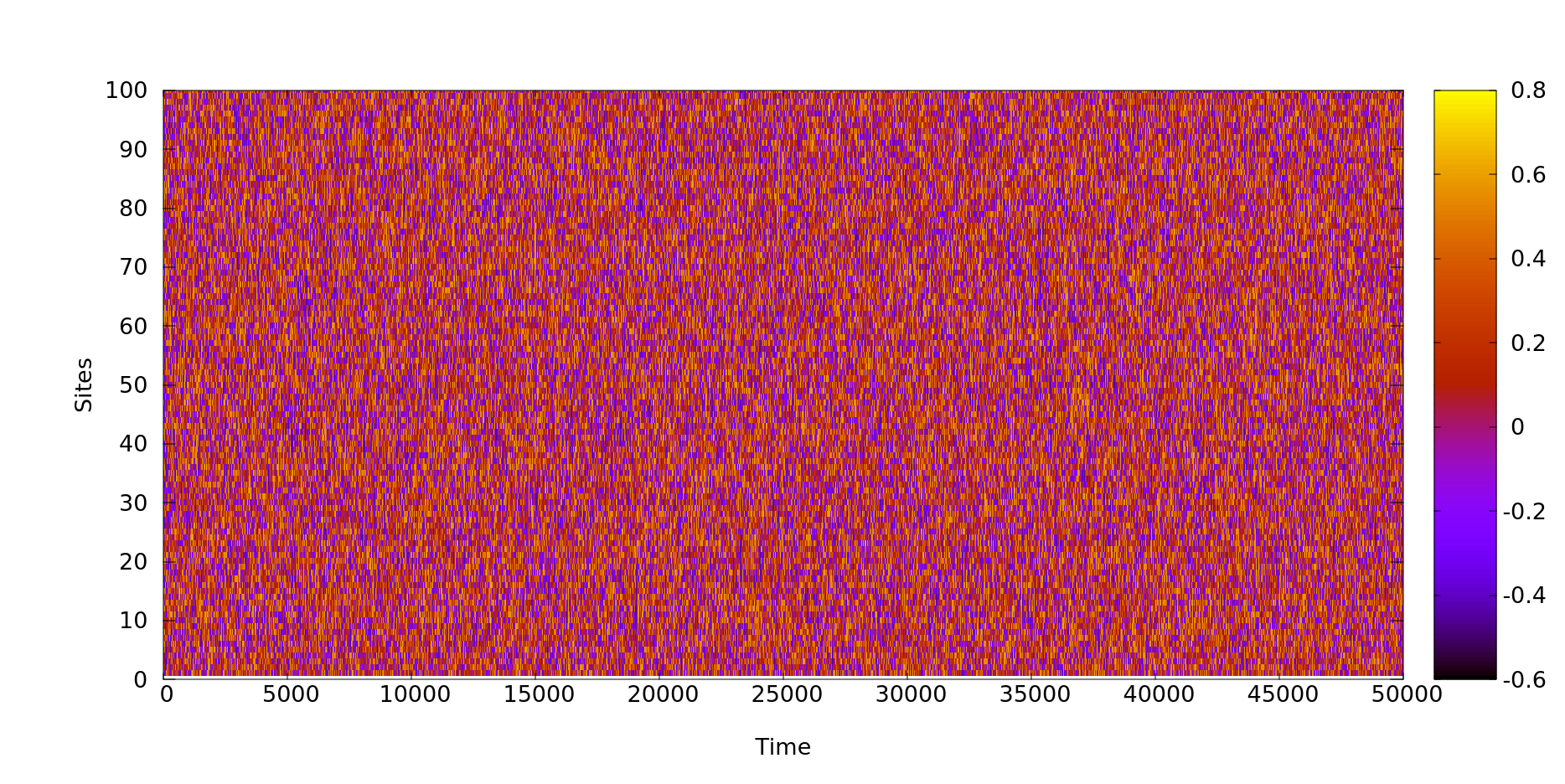} }}%
\qquad
\subfloat[]{{\includegraphics[width=\linewidth]{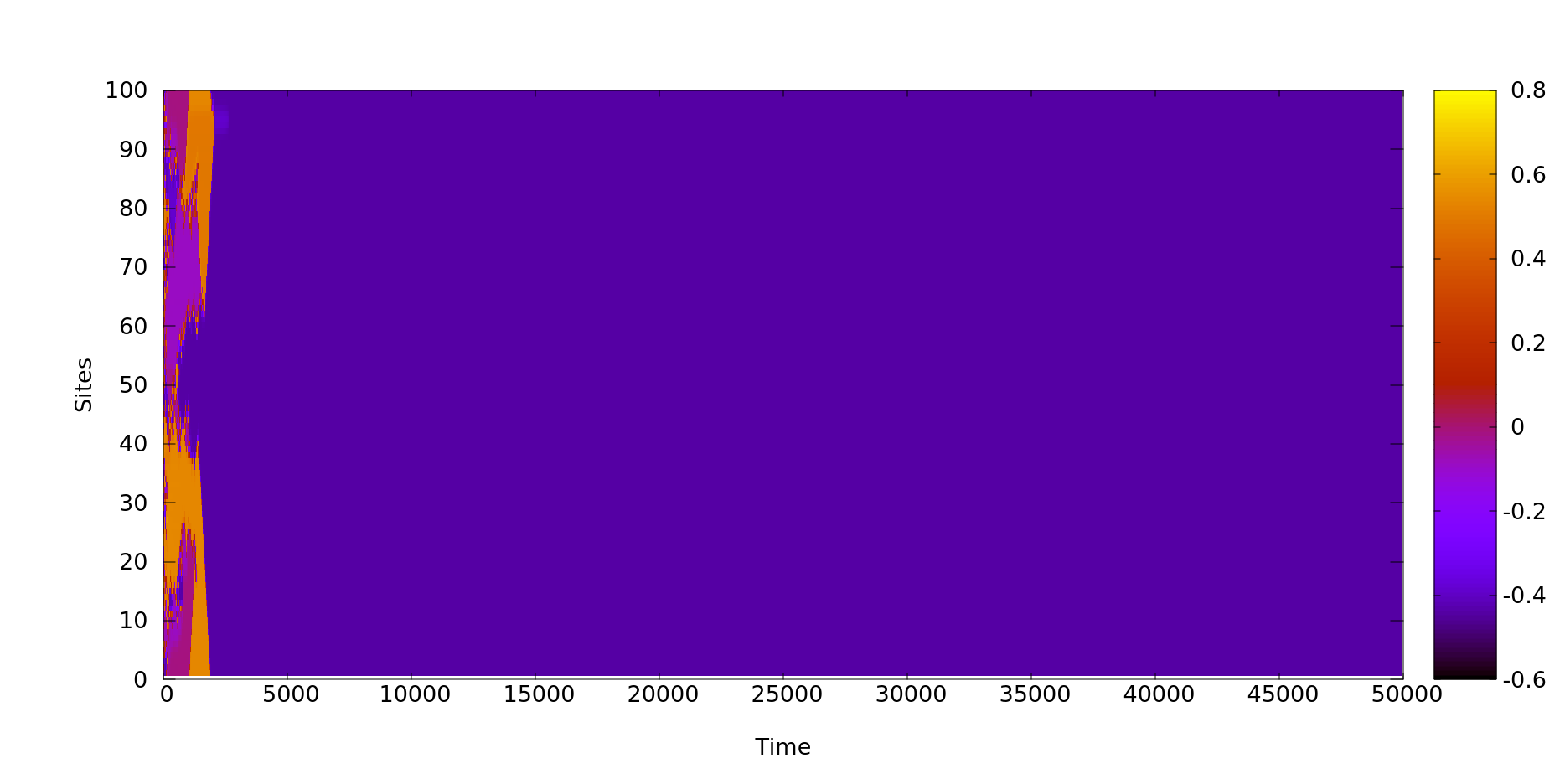} }}%
\caption{Heat map demonstrating spatiotemporal evolution of the system
for (a) $\alpha=0.8$ and $\beta=-0.40$
showing unsynchronized behavior for all time-steps while for  (b) $\alpha=0.8$ 
and $\beta=-0.45$ the system achieves synchronization early on.}%
\label{fig:3}%
\end{figure}

\begin{figure}%{.3\linewidth}
\centering
\subfloat[]{{\includegraphics[width=\linewidth]{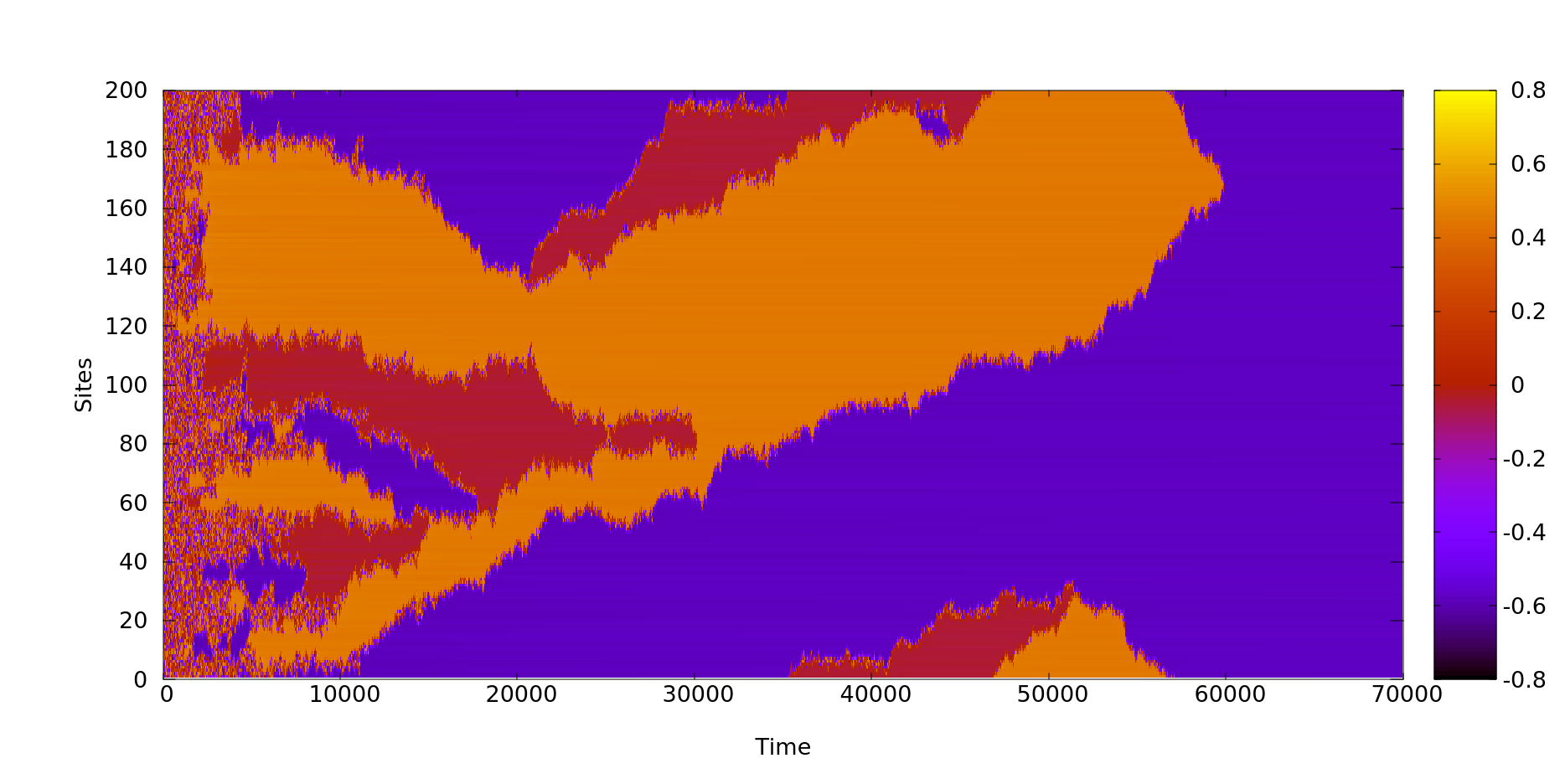} }}%
\qquad
\subfloat[]{{\includegraphics[width=\linewidth]{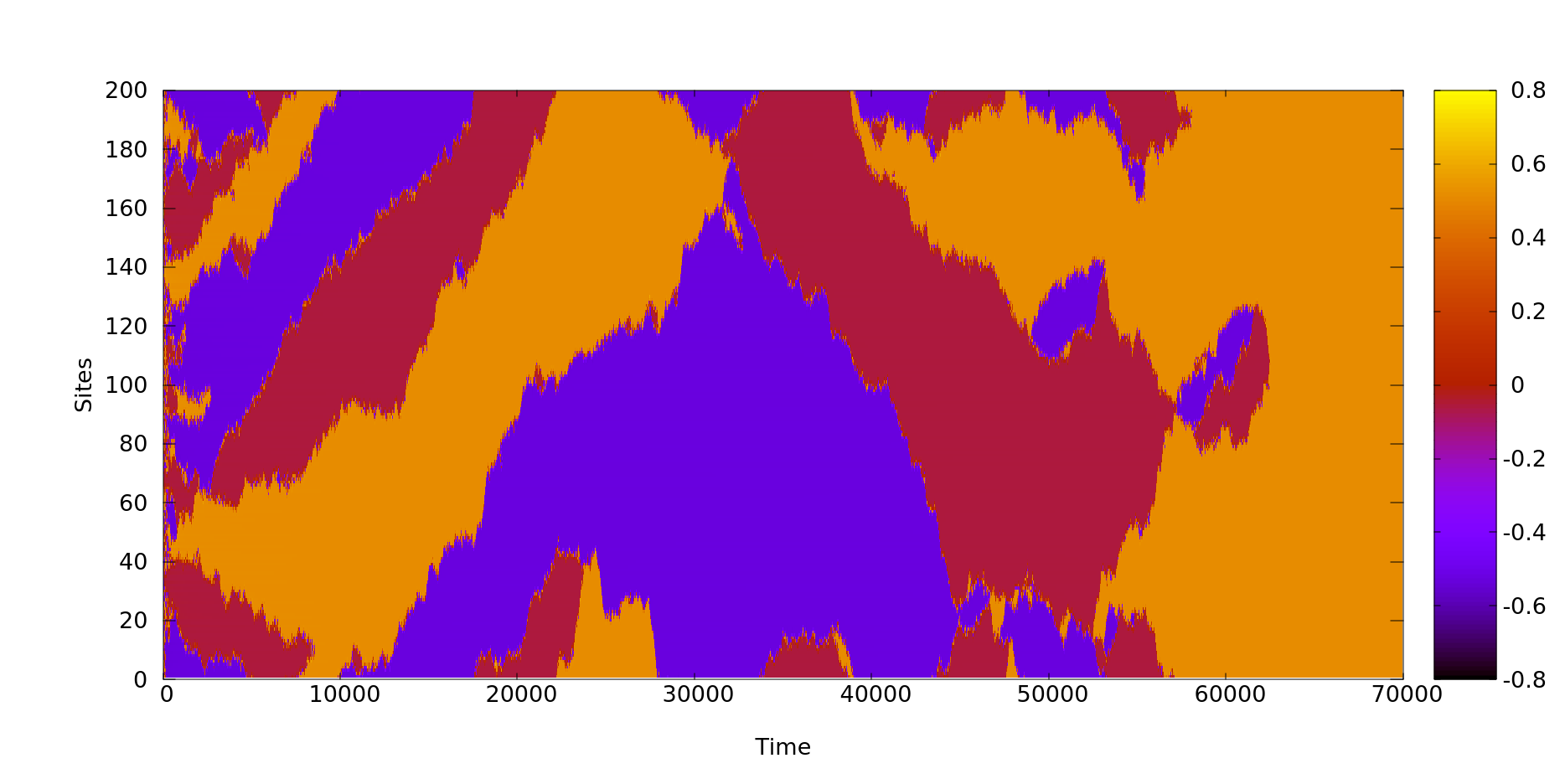} }}%
\caption{Heat map demonstrating spatiotemporal evolution of the system
	for (a) $\alpha=0.4$ and $\beta=-0.50$
	(b) $\alpha=0.6$ and $\beta=-0.44$. We simulate 200 maps and require $7 \times 10^4$
time-steps for achieving synchronization.}%
\label{fig:4}%
\end{figure}

\begin{figure}
\begin{center}
\includegraphics[width=85mm]{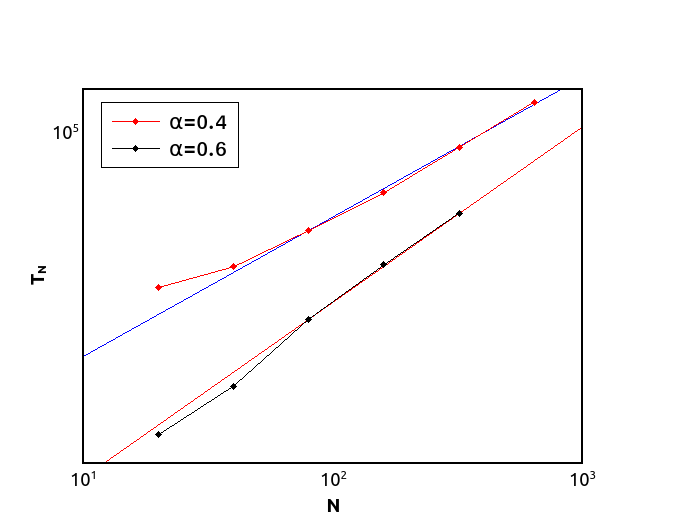}
\caption{Average time $T_N$ against system size $N$ is a  
power-law for $\alpha=0.4$ and $\alpha=0.6$.}%
\label{fig:5}
\end{center}
\end{figure}

Periodicity of synchronized patterns can be demonstrated by
plotting the time series of an arbitrary site in time.
The time series is shown in Fig. 6
 for the above values of $\alpha$. 
The dynamics is practically independent of initial conditions as
long as it is random.

\begin{figure}%
\centering
\subfloat[]{{\includegraphics[width=6cm]{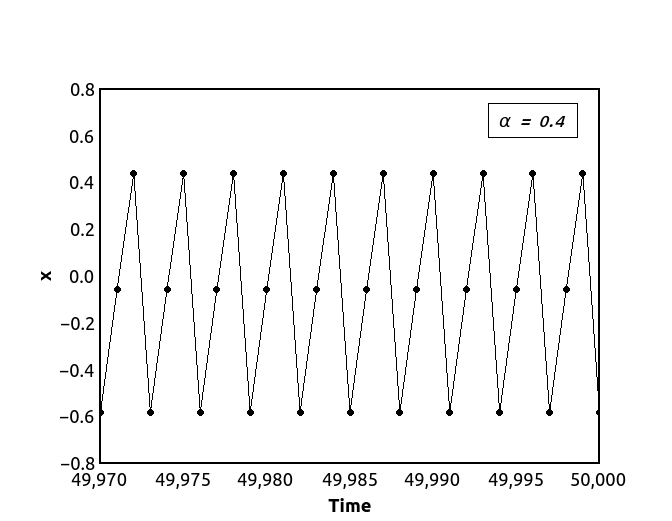} }}%
\qquad
\subfloat[]{{\includegraphics[width=6cm]{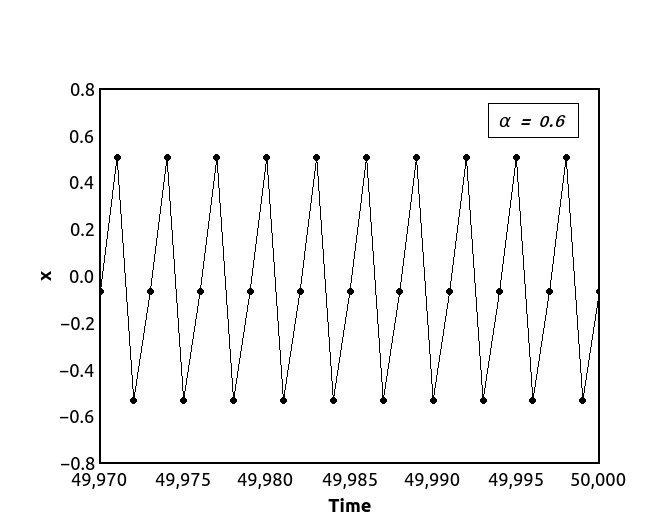} }}%
\qquad
\subfloat[]{{\includegraphics[width=6cm]{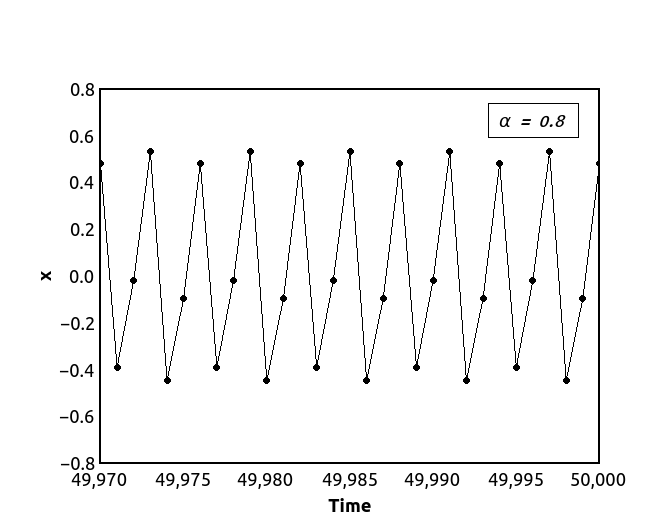} }}%
\qquad
	\caption{Synchronization with period-3 for (a) $\alpha=0.4$ and (b) $\alpha=0.6$
	(c) with period-6 for $\alpha=0.8$ }%
\label{fig:6}%
\end{figure}

\section{Globally Coupled Fractional Gauss Maps}

We can easily extend the above definition to the high-dimensional lattice. 
However, it would be cumbersome to simulate such a system and observe patterns in it.
Extension of this definition to  the mean-field model  may shed light
on behaviour in a very high dimensional system.
We define globally coupled fractional map as,

\begin{equation}
x(i,t) = x(i,0)+ \frac{1}{\varGamma(\alpha)} \sum_{j=1}^{t} g_{\alpha}(t-j) 
G_2\left(x(i,j-1), \sum_{k=1}^{N} f(x(k,j-1))\right),
\end{equation}
where,
\begin{equation}
G_2(a,b)=(1-\epsilon)f(a)+{\frac{\epsilon}{N}} b -a.
\end{equation}

We choose Gauss map as the map at a given site and vary the system parameters. 
We observe synchronization over a large range of parameters.
The standard deviation of spatial 
profile system decays as a power-law in time. It goes as
$t^{-\alpha}$ where $\alpha$ is the 
fractional-order parameter. For $\alpha= 0.4$ the power-law exponent is 0.4, 
for
$\alpha= 0.6$ it is 0.6 and for $\alpha= 0.8$ it is 0.8.
It is shown in Fig. 7.
This nature of decay of standard deviation depends only on fractional
order parameter and does not depend on other parameter values.
It
is superposed with oscillations  if the spatially synchronized state is
periodic in time.
For a period-2 synchronization,
the standard deviation
 takes smaller different values at odd and even
time-steps
for $\alpha=0.6$ and $\alpha=0.4$.
 These period-2 oscillations
 in the standard deviation 
for $\alpha=0.6$ are shown in Fig. 8(a).  In  Fig. 8(b),
we show a time series of some site. It clearly shows
2-period in time. Thus the period-2 oscillations in the standard deviation
are an artefact of period-2 in synchronized phase. 
If we ignore the oscillations,
the standard deviation decays with time as $t^{-\alpha}$.
This is an
expected behavior of function of Mittag-Leffler type
$e_{\alpha}(t)$ which is related to Mittag-Leffler function as
$e_{\alpha}(t):=E_{\alpha}(-t^{\alpha})$ for $t>0$. 
 This function  is important since 
it is an eigenfunction of the fractional relaxation equation.
Mainardi \cite{mainardi2013some} has shown that this function smoothly interpolates
between stretched exponential function at small times and
power-law decay at large times.
\begin{equation}
	E_{\alpha}(-t^{\alpha})=e_{\alpha}^{0}(t)=\exp\left[-\frac{t^{\alpha}}{\Gamma(1+\alpha)}\right], t\rightarrow0;
\end{equation}

\begin{equation}
	E_{\alpha}(-t^{\alpha})=e_{\alpha}^{\infty}(t)=\frac{t^{-\alpha}}{\Gamma(1-\alpha)}=\frac{\sin(\alpha\pi)}{\pi}\frac{\Gamma(\alpha)}{t^{\alpha}}, t\rightarrow\infty.
\end{equation}

\begin{figure}
\begin{center}
\includegraphics[width=85mm]{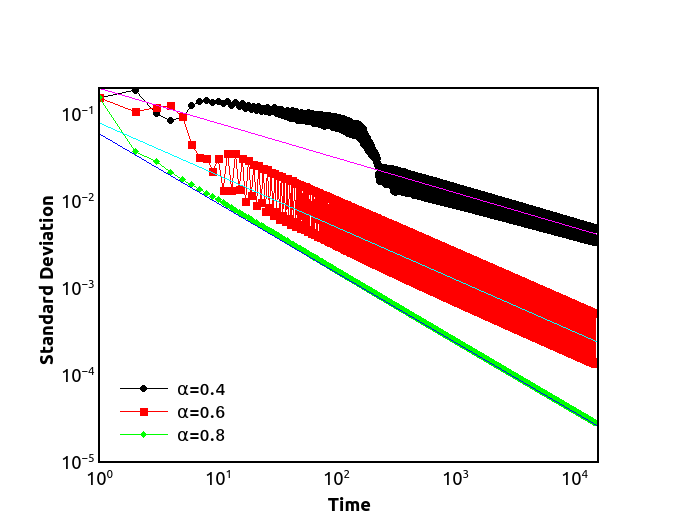}
\caption{Power-law decay of the standard deviation 
	in globally coupled fractional gauss 
	maps for $\epsilon=0.55$ and $\beta=-0.9$. 
	For $\alpha=0.4$ the decay exponent is $0.4$.
	For $\alpha=0.6$ the   exponent is $0.6$.
	For $\alpha=0.8$ the exponent is $0.8$.}
\label{fig:7}
\end{center}
\end{figure}

\begin{figure}
\centering
\subfloat[]{{\includegraphics[width=6cm]{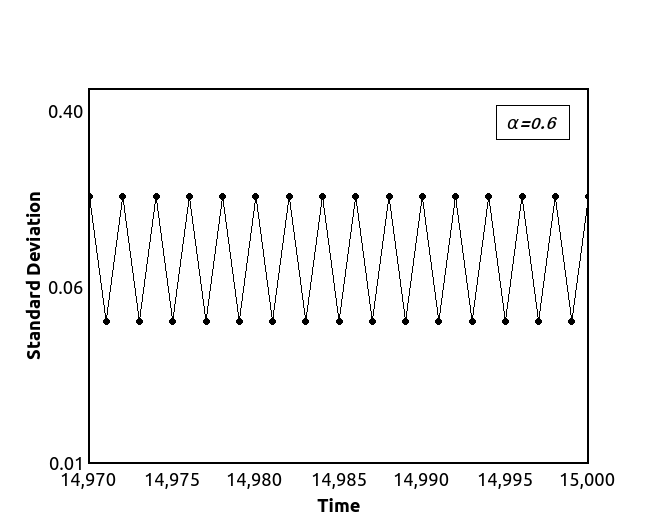} }}%
\qquad
\subfloat[]{{\includegraphics[width=6cm]{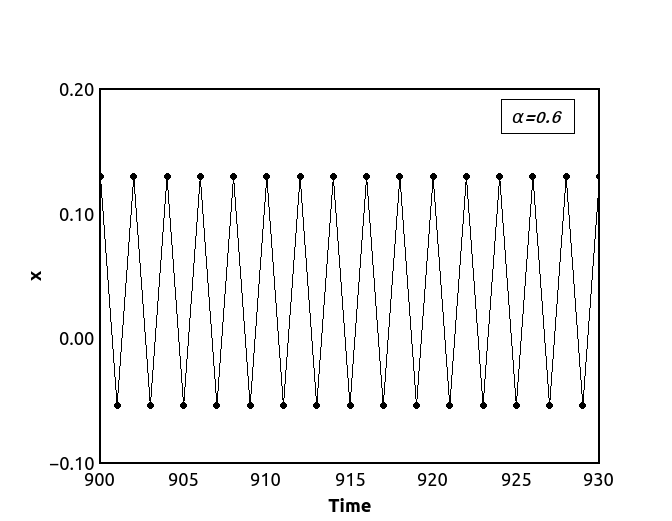} }}%
\qquad
\caption{Period-2 in globally coupled system for $\alpha=0.6$ in
	(a) standard deviation 
       (b) temporal evolution. }
\label{fig:8}
\end{figure}

\section{Coupled Fractional Gauss Maps with Small-World Networks}

Of late, dynamics on complex networks has been extensively investigated.
Most systems such as food web or internet connections do
not have a topology of d-dimensional Euclidean space. 
Small-world, as well as scale-free networks, have been very popular models of complex networks.
In small-world networks, we have a high clustering coefficient as well as
small path length. We study coupled fractional maps on small-world networks.
We begin with a one-dimensional chain where each element is coupled to
two nearest neighbours on either side and replace the nearest neighbour by a 
randomly chosen site on the lattice with probability $p$. For $p=1$,
it is a network with random (undirected) neighbours which is expected
to behave like a mean-field or globally coupled system. For $p=0$,
it is a simple one-dimensional chain. The small-world network allows us to interpolate between these two limits. For an equilibrium system,
any value of $p>0$ should have a behaviour similar to the mean-field model while
the situation is not so clear for nonequilibrium systems. 
We fix $p=0.7$ and
$\beta=-0.9$,
and investigate the behavior for various values of $\epsilon$ and $\alpha$.
The system is defined as follows, 

\begin{equation}
x(i,t) = x(i,0) + \frac{1}{\varGamma(\alpha)}\sum_{j=1}^{t}g\left(t-j\right)
G_3(x(i,j-1),
\sum_{k=1}^{4} f\left(x\left(\xi\left(
	i,k\right),j-1\right)\right))
,
\end{equation}
where, 
\begin{equation}
G_3(a,b) = (1-\epsilon)f(a)
+\frac{\epsilon}{4}
b -a.
\end{equation}
We consider a system with $N=4\times 10^3$ sites. We simulate this system for $10^4$ 
time steps for 
various values of $\alpha$ and $\epsilon$. 
We observe that the system gets synchronized for $\epsilon=0.8$.
The standard deviation of the 
system decays with time as
$t^{-\alpha}$ for all three values of $\alpha$.
This behaviour is similar to the decay in the globally coupled system
discussed in the above section. 
For $\alpha=0.6$ and $0.8$, and $\epsilon=1.0$,
same behavior is observed(See Fig. 9).

\begin{figure}
\begin{center}
\includegraphics[width=85mm]{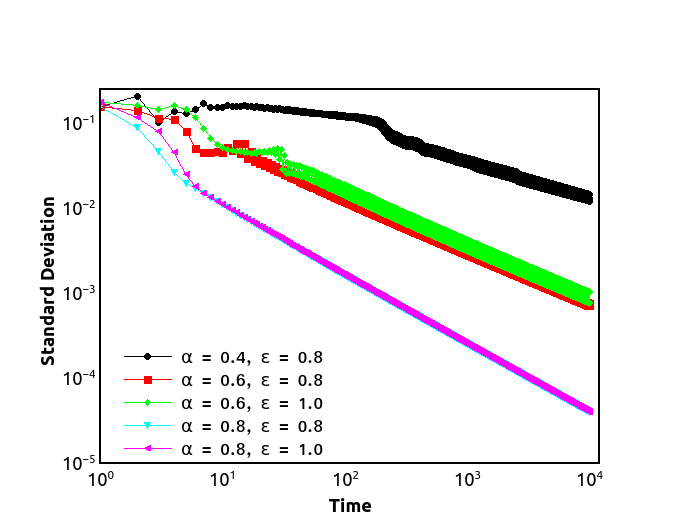}
\caption{Power-law decay of the standard deviation in the  system consisting of  small-world
 coupled fractional Gauss maps 	for various values of $\alpha$ and $\epsilon$.  
For $\alpha=0.4$ and $\epsilon=0.8$, the  exponent is $0.4$.
for $\alpha=0.6$ and $\epsilon=0.8$, the  exponent is $0.6$.
For $\alpha=0.6$ and $\epsilon=1.0$, the  exponent is $0.62$.
For $\alpha=0.8$ and $\epsilon=0.8$, the  exponent is $0.8$.
For $\alpha=0.8$ and $\epsilon=1.0$, the exponent  is $0.8$.}
\label{fig:9}
\end{center}
\end{figure}

\section{Results and Discussion}

We have extended the definition of fractional maps to a coupled map lattice
on arbitrary topology.
The system of coupled Gauss maps 
is analyzed for three different values of the fractional-order parameter, 
$\alpha$ and different values of map parameter $\beta$. 
We have presented results for locally coupled, globally coupled and small-world networks. 
For locally coupled Gauss maps, we observe spatially
synchronized states with period-3 or period-6 in time.  
Our studies indicate that such synchronization is possible
even in the thermodynamic limit.
For globally coupled Gauss maps, as well as for Gauss maps on a small-world network,
we observe synchronization. The approach to synchronization is non-exponential over a broad range of parameters. 
It is a clear power-law decay.
This approach is a lot slower than
what is observed in integer-order systems.
These results should be relevant in studies of 
spatially extended dynamical systems with
memory. In particular, it is of interest that a)
we can obtain synchronized periodic states 
in thermodynamic limit even
with local coupling and b)
approach to synchronization is a power-law with
power same as the fractional order parameter. 
We note that for integer-order coupled map lattice,
the standard deviation decays exponentially in the synchronized phase.
and
power-law decay is obtained only at the critical point for continuous
transitions.

For globally coupled maps as well as for
the small-world  topology of coupled map lattices, synchronization is
obtained for larger parameter values.
Power-laws in space and time are observed in a range of
physical phenomena and are not so common in
mathematical models. A popular model for power-laws in space and
time has been self-organized criticality \cite{bak2013nature}. 
Of late, the Griffiths phase \cite{moretti2013griffiths}
has been offered as an explanation for power-laws observed in
systems such as coupled neurons. Above studies point to another
possibility, namely, memory in the physical system leads to
long-time dynamical correlations.

\section{Acknowledgments}
PMG thanks DST for financial help (EMR/2016/006685). 
 \bibliographystyle{elsarticle-num} 
 \bibliography{edited.bib}

\end{document}